\documentclass[12pt]{amsart}
\usepackage[latin1]{inputenc}
\usepackage{amsmath,latexsym,amssymb,graphicx,dsfont}
\usepackage{natbib}

\newcommand{\R}{\mathbb{R}}

\newtheorem{Th}{Theorem}[section]
\newtheorem{Lemma}[Th]{Lemma}
\newtheorem{Cor}[Th]{Corollary}
\newtheorem{Prop}[Th]{Proposition}

\newcommand{\E}{\mathbb{E}}
\newcommand{\Prob}{\mathbb{P}}
\newcommand{\calL}{\mathcal{L}}

\newcommand{\calN}{\mathcal{N}}

\newcommand{\calR}{\mathcal{R}}
\newcommand{\N}{\mathbb{N}}

\newcommand{\Z}{\mathbb{Z}}

\renewcommand{\H}[1]{\mathbf{H}_{#1}}
\newcommand{\s}[1]{\mathbf{s}_{#1}}
\renewcommand{\t}[1]{\mathbf{t}_{#1}}

\numberwithin{equation}{section}

\begin{document}

\title{Acceleration of Lamplighter Random Walks}

\address{Institut f\"ur mathematische Strukturtheorie (Math. C), Graz University of Technology, Steyrergasse
  30, A-8010 Graz, Austria}

\email{gilch@TUGraz.at}
\urladdr{http://www.math.tugraz.at/$\sim$gilch/}
\date{\today}

\maketitle

\centerline{\scshape Lorenz A. Gilch}
\medskip
{\footnotesize
 \centerline{Graz University of Technology, Graz, Austria}}

\begin{abstract}
Suppose we are given an infinite, finitely generated group $G$ and a transient
random walk on the wreath product $(\Z/ 2\Z)\wr G$, such that its projection on
$G$ is transient and has finite first moment. This random walk can be
interpreted as a lamplighter random walk on $G$. 
Our aim is to show that the random walk on the wreath
product escapes to infinity with respect to a suitable (pseudo-)metric faster 
than its projection onto $G$. We also address the case where the pseudo-metric is the length of a shortest ``travelling salesman
tour''. In this context, and excluding some degenerate cases if $G=\Z$,
the linear rate of escape is strictly  bigger than the rate of escape of
the lamplighter random walk's projection on $G$.
\end{abstract}

\section{Introduction}


Let $G$ be an infinite group generated by a finite symmetric set $S$, and
imagine a lamp sitting at each group element. These lamps have two states: 0
(``off'') or 1 (``on''), and initially all lamps are off. We think of a
lamplighter walking randomly on $G$ and switching lamps on or off as he walks. We investigate the
following model: at each step
the lamplighter may walk to some random neighbour vertex, and may flip some lamps in a
bounded neighbourhood of his position. This model can be interpreted as
a random walk on the wreath product $(\Z/2\Z)\wr G$ governed by a probability
measure $\mu$. The random walk
is described by a Markov chain $(Z_n)_{n\in\N_0}$, which represents the position $X_n$ of the
 lamplighter and the lamp configuration $\eta_n:G\to \Z/2\Z$ at time $n$. We assume that the
 lamplighter random walk's projection on $G$ has finite first moment and is also transient.
\par
For better visualization, we identify $G$ with its Cayley graph with
respect to the generating set $S$.
Suppose we are given ``lengths'' of the elements of $S$ such that $s\in S$ and
$s^{-1}\in S$ have the same length. 
The length of a path in
 $G$ is the sum of the lengths of its edges, and we denote by
 $d(\cdot,\cdot)$ the metric on $G$ induced by the lengths of the edges. 
We denote by $d_\mathrm{TS}(\eta,x)$ the length of an optimal ``travelling
salesman tour'' from the identity $e$ to $x\in G$ that visits each point in
$\mathrm{supp}(\eta)$ (where $\eta:G\to \Z/2\Z$ has finite support).
Then a natural length function
 $\ell(\eta,x)$ is given by $d_\mathrm{TS}(\eta,x)+c_\calL \cdot
 |\mathrm{supp}(\eta)|$ for an arbitrary, but fixed constant $c_\calL\geq
 0$. By transience, our random walk escapes to
 infinity with respect to this length function. 
\par
The (new) topic that we address in this paper is the
comparison of the limits
 $\ell=\lim_{n\to\infty}\ell(Z_n)/n$ and $\ell_0=\lim_{n\to\infty} d(e,X_n)/n$,
which exist almost surely. They describe the speed of the lamplighter random walk and its projection on
 $G$, respectively.
The number $\ell_0$ is called the
 \textit{rate of escape}, or the \textit{drift}, of $(X_n)_{n\in\N_0}$ and $\ell$ is the
 \textit{rate of escape of the lamplighter random walk} $(Z_n)_{n\in\N_0}$. 
It is well-known that the rate of escape exists for random walks with finite
first moment on transitive graphs. This follows from \textit{Kingman's subadditive
  ergodic  theorem}; see Kingman \cite{kingman}, Derriennic \cite{derrienic} and
\mbox{Guivarc'h \cite{guivarch}.} 
We will prove that, under some weak assumptions on
$G$, we have $\ell>\ell_0$, that is, the
lamplighter random walk escapes strictly faster
to infinity than its projection onto $G$, on which we have the
 metric $d(\cdot,\cdot)$. If the lamplighter random walk's projection on $G$
is transient and has zero drift, then the acceleration of the lamplighter
random walk follows from results of Kaimanovich and Vershik \cite{kaimanovich-vershik} and of
Varopoulos \cite{varopoulos}. Thus, we may restrict ourselves to the case $\ell_0>0$. More explicitly, we will prove that
 $\lim_{n\to\infty} |\mathrm{supp}(\eta_n)|/n>0$, where $\eta_n$ is the lamp
 configuration at time $n$. From this it follows directly that $\ell>\ell_0$. We also
 prove \mbox{$\lim_{n\to\infty} d_\mathrm{TS}(Z_n)/n>\ell_0$} (except for some
 degenerate cases), providing $\ell>\ell_0$. 
\par
Let us briefly review a few selected results regarding the rate of escape.
The classical case is that of random walks on the $k$-dimensional grid $\mathbb{Z}^k$, where $k\geq 1$, which can be described by the sum of $n$
i.i.d. random variables, the increments of $n$ steps. By the law of
large numbers the limit $\lim_{n\to\infty} \Vert Z_n \Vert /n$, where $\Vert \cdot\Vert$ is the
distance on the grid to the starting point of the random walk, exists almost
surely. Furthermore, this limit is positive if the increments have non-zero
mean vector. There is an important link
between drift and the Liouville property: the entropy (introduced by Avez \cite{avez72}) of any random walk on a
group is non-zero if and only if non-constant harmonic functions exist; see
Kaimanovich and Vershik \cite{kaimanovich-vershik} and Derriennic
\cite{derrienic}. Moreover, if the rate of escape is zero, then the entropy is
zero (first observed by Guivarc'h \cite{guivarch}). Varopoulos
\cite{varopoulos} has shown the converse for symmetric finite range random
walks on groups. The recent work of Karlsson and
Ledrappier \cite{karlsson-ledrappier07} generalizes this result to symmetric
random walks with finite first moment of the step lengths.
\par
In this paper we deal with random walks on wreath products, for which there
are many detailed results: Lyons,
Pemantle and Peres \cite{lyons-pemantle-peres} gave a lower bound for the rate
of escape of inward-biased random walks on lamplighter groups. Revelle \cite{revelle} examined the rate of escape of random walks on wreath
products. He proved laws of the iterated logarithm for the inner and outer
radius of escape.  For a finitely generated group $A$, Dyubina
\cite{dyubina} proved that the drift w.r.t. the word metric of a random walk on the wreath product $(\mathbb{Z}/2\Z)\wr
A$ is
zero if and only if the random walk's projection onto $A$ is recurrent.   
\par
It is not obvious that lamplighter random walks are in general faster than
their projections onto $G$: e.g., consider the Switch-Walk lamplighter random walk on
$\mathbb{Z}$ with drift: in each step switch the lamp at the actual position with
probability $p\in(0,1)$ and then walk to a random neighbour vertex. Then the rate of escape of the lamplighter random walk
is equal to the random walk's projection onto $\mathbb{Z}$ whenever
$c_\calL=0$; compare with Bertacchi \cite{bertacchi}.  However, this example is
more or less the only counterexample. The author of this article has investigated the rate of escape
of lamplighter random walks arising from a simple random walk on homogeneous trees providing tight lower and
upper bounds for the rate of escape; see Gilch \cite{gilch2}. In particular,
the lamplighter random walk is significantly faster than its projection onto the tree. This was the
starting point for the investigation of the relation between $\ell$ and
$\ell_0$ on more general classes of graphs.
\par
The structure of this article is as follows: in Section \ref{lamplighter-graphs} we give an
introduction to lamplighter random walks on groups and some basic
properties. 
In Section \ref{WoS} we prove that $\lim_{n\to\infty}
|\mathrm{supp}(\eta_n)|/n>0$ if and only if the lamplighter random walk's projection on $G$
is transient. In Section \ref{SWS} we
prove $\lim_{n\to\infty} d_\mathrm{TS}(Z_n)/n>\ell_0$ under some weak
assumptions on $S$, which exclude some degenerate cases. 
Finally, in Section \ref{remarks} we give some
additional remarks regarding extensions of the presented results.

\section{Lamplighter Groups}

\label{lamplighter-graphs}
\subsection{Groups and Random Walks}

Consider an infinite, finitely generated group $G$ with identity $e$ and a
finite, symmetric set of generators $S\subseteq G\setminus \{e\}$, which
generates $G$ as a semigroup. We assign to each $s\in S$ a \textit{length} $l(s)=l(s{^{-1}})>0$. We write $r_1:=\min_{s\in S} l(s)$. These lengths induce a metric on $G$: the distance between $x,y\in G$ is given by
$$
d(x,y) := \min \biggl\lbrace
\sum_{i=1}^n l(s_i)\, \Bigl| \, s_1,\dots,s_n\in S \textrm{ such that } y=xs_1s_2\cdots s_n
\biggr\rbrace.
$$
We identify $G$ with its Cayley graph with respect to $S$.
A \textit{path} in $G$ is a finite sequence of group elements $[x_0,x_1,\dots,x_n]$ such that $x_{i-1}^{-1}x_i\in S$. The \textit{length} of this path is $\sum_{i=1}^n l(x_{i-1}^{-1}x_i)$.
The \textit{ball} $B(x,r)$ centered at $x\in G$ with radius
$r\geq 0$ is given
by the set of all elements $y\in G$ with $d(x,y)\leq r$.

\subsection{Lamplighter Random Walks}
\label{LRW}
Imagine a \textit{lamp} sitting at each vertex of $G$, which can be switched
off or on, encoded by ``0'' and ``1''. We think of a lamplighter walking on
$G$ and switching lamps on and off. The lamp configurations are encoded by functions $\eta:G\to\Z/2\Z$. Writing $\Z_2:=\Z/2\Z$, the \textit{set of finitely supported
configurations} of lamps is
$$
\calN := \bigl\lbrace \eta: G\to \mathbb{Z}_2\ \bigl| \ |\mathrm{supp}(\eta)|<\infty\bigr\rbrace.
$$
Here, $|A|$ denotes the cardinality of a set $A$. Denote by $\mathbf{0}$ the
zero configuration, which will be the initial lamp configuration of the random walk, and by $\mathds{1}_x$ the
configuration where only the lamp at $x\in G$ is on and all other lamps are off. 
The \textit{wreath product} of $\mathbb{Z}_2$ with $G$ is
$$
\mathcal{L} :=  \Bigl( \sum_{x\in G} \mathbb{Z}_2 \Bigr)  \rtimes G = \mathbb{Z}_2 \wr G. 
$$
The elements of $\calL$ are pairs of the form $(\eta,x)\in \calN\times
G$, where $\eta$ represents a configuration of the lamps and $x$ the position
of the lamplighter. For \mbox{$x,w\in G$} and $\eta\in\calN$, define
$$
(x\eta)(w):=\eta(x^{-1}w).
$$
A group operation on $\mathcal{L}$ is given by
$$
(\eta_1,x)(\eta_2,y) := \bigl( \eta_1 \oplus (x\eta_2),xy\bigr),
$$
where $x,y\in G$, $\eta_1,\eta_2\in \calN$, $\oplus$ is
the componentwise addition modulo 2. The group identity is $(\mathbf{0},e)$. We
call $\calL$ together with this operation the \textit{lamplighter group} \mbox{\textit{over} $G$.}
\par
A natural symmetric set of generators of $\calL$ is given by
$$
S_{\mathcal{L}}:=\bigl\lbrace (\mathds{1}_e,e),(\mathbf{0},s)
\, \bigl|\, s\in S \bigr\rbrace.
$$
Consider the Cayley graph of $\mathcal{L}$ with respect to $S_{\mathcal{L}}$. We lift
$d(\cdot,\cdot)$ to a \mbox{(pseudo-)metric} $d_{\calL}(\cdot,\cdot)$ on $\calL$ by
assigning the following lengths to the elements of $S_\calL$: 
$l\bigl((\mathbf{0},s)\bigr):=l(s)$ for $s\in S$ and $l\bigl((\mathds{1}_e,e)\bigr):=c_\calL\geq 0$, where $c_\calL$ is some arbitrary, but fixed non-negative constant. These lengths induce a (pseudo-)metric on $\calL$.
 The distance between
$(\eta,x)$ and $(\eta',y)$ is then the minimal length of all paths in the
Cayley graph of $\calL$ joining both vertices. More explicitely, we denote by
$d_{\mathrm{TS}}(\eta,x)$ the minimal length of a ``travelling salesman tour'' on $G$
(not on $\calL$) from $e$ to $x$, which visits each \mbox{$y\in\mathrm{supp}(\eta)$.} With this notation, we have 
$$
\ell(\eta,x) :=  d_\calL\bigl((\mathbf{0},e),(\eta,x)\bigr) = d_{\mathrm{TS}}(\eta,x) + c_\calL \cdot |\mathrm{supp}(\eta)|.
$$
The case $c_\calL=0$ can also be interpreted as the model where $S_\calL$ is
replaced by $\{(\mathbf{0},s),(\mathds{1}_e,s)\mid s\in S\}$ and where the
length of $(\mathds{1}_e,s)$ equals $l(s)$. In this case, lamp switches are not
charged by the pseudo-metric.
\par
We now consider an irreducible, transient random walk on $\calL$ starting at the identity $(\mathbf{0},e)$ such that the random walk's
projection onto $G$ is also transient. For this purpose, consider the
sequence of i.i.d. $\calL$-valued random variables
$(\mathbf{i}_n)_{n\in\mathbb{N}}$ governed by a probability
measure $\mu$ which satisfies the following conditions:
\begin{enumerate}
\item[(1)] $\langle \mathrm{supp}(\mu)\rangle =\calL$.
\item[(2)] There is a non-negative real number $R$ such that 
$$
\mu(\eta,x)>0\ \textrm{ implies } \ d(e,y)\leq R  \textrm{ for all } y\in \mathrm{supp}(\eta).
$$
\item[(3)] The projection of $\mu$ onto $G$ has finite first moment, that is,\\
  \mbox{$\sum_{(\eta,x)\in\calL} d(e,x)\mu(\eta,x)<\infty$.}
\end{enumerate}
We write $\mu^{(n)}$ for the $n$-th convolution power of $\mu$.
A lamplighter
random walk starting at $(\mathbf{0},e)$ is described by the sequence of $\calL$-valued random variables $(Z_n)_{n\in\mathbb{N}_0}$ in the
following natural way:
$$
Z_0:=(\mathbf{0},e), \quad Z_n:=Z_{n-1}\mathbf{i}_n \ \textrm{ for all } n\geq 1.
$$
More precisely, we write \mbox{$Z_n=(\eta_n,X_n)$}, where $\eta_n$ is the random configuration of the lamps at
time $n$ and $X_n$ is the random group element at which the lamplighter
stands  at time $n$. As a \textit{general assumption} we assume transience of
$(X_n)_{n\in\mathbb{N}_0}$. We explain below what happens if this assumption fails.
\par
The corresponding single and $n$-step transition probabilities of the random
walk on $\calL$ are denoted by $p(\cdot,\cdot)$ and
$p^{(n)}(\cdot,\cdot)$. 
We write $\Prob_z[\,\cdot\,]:=\Prob[\,\cdot \mid Z_0=z]$
for any $z\in\calL$, if we want to start the lamplighter walk at $z$ instead of
$(\mathbf{0},e)$.
\par
Observe that by transience of
$(X_n)_{n\in\mathbb{N}_0}$ each finite subset of $G$ is visited only finitely often yielding
that the sequence $(\eta_n)_{n\in\N_0}$ converges pointwise to a random limit configuration
\mbox{$\eta_\infty: G\to\mathbb{Z}_2$}, which is not necessarily finitely
supported. On the other hand, $(X_n)_{n\in\N_0}$ leaves every finite set after some finite time forever, that is, $d(e,X_n)$ goes to infinity.
\par
As a consequence of \textit{Kingman's subadditive ergodic theorem} there are
non-negative numbers \mbox{$\ell_0,\ell \in\R$} such that
$$
\ell_0=\lim_{n\to\infty} \frac{d(e,X_n)}{n}\ \textrm{ and } \
\ell=\lim_{n\to\infty} \frac{\ell(Z_n)}{n}\ 
\textrm{ almost surely;}
$$
see Derriennic \cite{derrienic} and Guivarc'h \cite{guivarch}. The number $\ell_0$ is  called the \textit{rate of escape} or \textit{drift} of the lamplighter random
walk's projection onto $G$. Analogously, $\ell$ is the \textit{rate of escape
  of the lamplighter random walk}. Moreover, we can write 
\begin{eqnarray}\label{def-ell}
\ell = \ell_{\mathrm{TS}} +c_\calL \cdot \ell_{\mathrm{supp}},
\end{eqnarray}
where 
\begin{eqnarray*}
\ell_{\mathrm{TS}}&:= &\lim_{n\to\infty} \frac{d_{\mathrm{TS}}(\eta_n,X_n)}{n} \ \textrm{ (core rate
  of escape) and}\\
\ell_{\mathrm{supp}}&:= &\lim_{n\to\infty} \frac{|\textrm{supp}(\eta)|}{n} \
\textrm{(asymptotic configuration size)}.
\end{eqnarray*}
The latter limits exist for the same reason as above. Obviously, $\ell\geq \ell_{\mathrm{TS}}\geq
\ell_0$. By Kaimanovich and Vershik
\cite{kaimanovich-vershik}, we have $h\leq \ell \cdot g$, where
$$
h=\lim_{n\to\infty} -\frac{1}{n} \sum_{(\eta,x)\in \calN\times G}
p^{(n)}\bigl((\mathbf{0},e),(\eta,x)\bigr)\, \log
p^{(n)}\bigl((\mathbf{0},e),(\eta,x)\bigr)
$$
is the \textit{asymptotic entropy} and
$g=\lim_{n\to\infty} (\log |B_\calL (n)|)/n$  is the \textit{growth rate} of
$\calL$, where $B_\calL(n)$ denotes the ball around $(\mathbf{0},e)$ of
radius $n$ with respect to $d_\calL(\cdot,\cdot)$. Existence of $h$ and $g$ follows again from Kingman's subadditive
ergodic theorem. We have $g<\infty$ even if one considers the balls with respect
to the pseudo-metric $d_\mathrm{TS}(\cdot, \cdot)$, because of the relation
$$
\bigl\lbrace (\eta,x)\in\calL \mid d_{\mathrm{TS}}(\eta,x)\leq n\bigr\rbrace \subseteq
B_\calL \bigl(n+(\lfloor n/r_1\rfloor +1) c_\calL\bigr).
$$
Furthermore, we also have $h>0$: the mapping
$$
(\eta,x) \mapsto \Prob_{(\eta,x)}[\eta_\infty(e)=0]
$$
defines a non-constant bounded harmonic function. Thus, the Poisson boundary is
non-trivial, that is, $h>0$; see Kaimanovich \cite{kaimanovich} 
and Kaimanovich and Woess \cite{kaimanovich-woess}.
Thus, we get $\ell>0$. As an additional remark, let us mention that Dyubina
\cite{dyubina} proved that lamplighter random walks w.r.t. the word metric have
non-zero drift
 if
and only if the projection on $G$ is transient. Since $S$ is finite and the
lengths of edges are bounded, we have in our situation $\ell=0$ if $G$ is recurrent.
\par
Our basic aim is to show
that $\ell$ is strictly bigger than $\ell_0$, that is, the lamplighter random walk
escapes faster to infinity than its projection onto $G$. For this purpose, we
will show $\ell_{\mathrm{supp}}>0$ in the following section, giving $\ell>\ell_0$ in the case
$c_\calL>0$. In Section \ref{SWS}, we will prove $\ell_{\mathrm{TS}}>\ell_0$
under suitable weak assumptions on $G$, which exclude degenerate cases for
special choices of $S$ and $l(\cdot)$ when $G=\Z$.

\section{The Asymptotic Configuration Size}

\label{WoS}

In this section we want to show that the number of lamps which are on increases
asymptotically at
linear speed.
 We show that $\ell_{\mathrm{supp}}>0$, giving $\ell>\ell_0$ in the case
 $c_\calL>0$. 
\par
Consider now the lamplighter random walk's projection on $G$ and its range
$\mathcal{R}_n\subseteq G$, which is the set of visited elements up to time
$n$. By Derriennic \cite{derrienic}, $|\mathcal{R}_n|/n$ converges to $\Prob[\forall n\geq 1:
X_n\neq e]$, which is strictly positive in our case by transience. For
$j\in\N$, let  
$$
\s{j} := \min \bigl\lbrace n\in\N_0 \,\bigl|\, |\mathcal{R}_n| = j \bigr\rbrace. 
$$
By transience, $\s{j}<\infty$ almost surely. For
$k\in\N_0$, let 
$$
\Delta_{k,j} := 
\begin{cases}
1, & \textrm{ if } \eta_{\s{j}+k}(X_{\s{j}})=1, \\
0, & \textrm{ otherwise}.
\end{cases}
$$
With this definition we have for all $n\in\N_0$
\begin{equation}\label{abschaetzung1}
|\textrm{supp}(\eta_{n})|
\geq \sum_{j=1}^{|\mathcal{R}_n|} \Delta_{n-\s{j},j}.
\end{equation}
We give now a uniform lower bound for the probability
$\Prob_{(\eta,x)}[\eta_\infty(x)=1]$ with \mbox{$(\eta,x)\in\calN\times G$}.
\begin{Lemma}\label{q-lemma}
There are $\kappa\in\N$ and $C>0$ such that 
for all $(\eta,x)\in\calN\times G$
$$
\Prob_{(\eta,x)}[\forall n\geq \kappa: \eta_n(x)=1] \geq C.
$$
\end{Lemma}
\begin{proof}
By transience and bounded range of the random walk $(X_n)_{n\in\N_0}$, there is at least one vertex $y\in G$ with \mbox{$R<d(e,y)$} such that 
$$
\tilde p:=\Prob_{(\eta',y)}[\forall n\geq 1: X_n\notin B(e,R) ]>0 \quad
\textrm{ for each } \eta'\in\calN.
$$
Moreover, there are $\kappa_0,\kappa_1\in\N$ such that
$C_0:=\mu^{(\kappa_0)}\bigl((\mathbf{0},y)\bigr)>0$ and $C_1:=\mu^{(\kappa_1)}\bigl((\mathds{1}_e,y)\bigr)>0$. Transitivity provides
that the probability of walking from $(\eta,x)$ to some
$(\eta',xy)$ with \mbox{$x\in\mathrm{supp}(\eta')$} in at most
$\kappa:=\max\{\kappa_0,\kappa_1\}$ steps is at least
$C':=\min\{C_0, C_1\}$. With $C:=C'\cdot \tilde p$ follows the claim of the lemma. 
\end{proof}
The next lemma gives a non-trivial uniform lower bound for $\E[\Delta_{n,j}]$
with $n\geq \kappa$:
\begin{Lemma}\label{delta-lemma1} 
For all $j,n\in \N$ with $n\geq \kappa$ we have $\E[\Delta_{n,j}]\geq C>0$.
\end{Lemma}
\begin{proof}
In order to bound $\Prob[\Delta_{n,j}=1]$ uniformly from below, we decompose
according to all possible states of $X_{\s{j}}$, followed by walking steps to some 
$X_{\s{j}}y\in G\setminus B(X_{\s{j}},R)$, such that the lamp at $X_{\s{j}}$ --
if necessary -- will be switched on and after reaching $X_{\s{j}}y$ the random walk
does not return to $B(X_{\s{j}},R)$.
\par
By vertex-transitivity and Lemma \ref{q-lemma}, the probability of starting
in $X_{\s{j}}$  and walking in at most $\kappa$ steps to
some vertex $X_{\s{j}}y\in G\setminus B(X_{\s{j}},R)$, such that the lamp at
$X_{\s{j}}$ rests on and no further visit in $B(X_{\s{j}},R)$
after reaching $X_{\s{j}}y$ is at least $C>0$.
\par
Observe that, by transience we have for each $j\in\N$
$$
\sum_{x\in G} \sum_{m\geq 0}
\Prob[\s{j}=m, X_m=x]=\Prob[ \s{j}<\infty]=1. 
$$
We get for $n\geq \kappa$:
\begin{eqnarray*}
\E[\Delta_{n,j}]
& = & \sum_{(\eta,x)\in \calL} \sum_{m\geq 0} 
\Prob\bigl[
\s{j}=m,Z_m=(\eta,x),\eta_{m+n}(x)=1\bigr]\\
&\geq & \sum_{(\eta,x)\in \calL} \sum_{m\geq 0} 
\Prob\bigl[\s{j}=m,Z_m=(\eta,x)\bigr] \cdot 
\Prob_{(\eta,x)}\bigl[\eta_{n}(x)=1\bigr]\\
&\geq & \Prob[ \s{j}<\infty] \cdot C = C >0.
\end{eqnarray*}
\end{proof}
Now we can conclude:
\begin{Th}\label{WoS-Th}
For the lamplighter random walk with respect to the infinite, finitely generated group $G$, 
$$
\ell_{\mathrm{supp}} >0.
$$
Moreover, if $c_\calL>0$, then $\ell>\ell_0$.
\end{Th}
\begin{proof}
Recall the definition of
$\kappa$ from Lemma \ref{q-lemma}. We have:
$$
\ell_{\mathrm{supp}} = \lim_{n\to\infty}  
\frac{|\mathrm{supp}(\eta_{n+\kappa})|}{|\calR_n|}\frac{|\calR_n|}{n}
\frac{n}{n+\kappa}
= \ell_1 \cdot \bar F,
$$
where $\ell_1:=\lim_{n\to\infty} |\mathrm{supp}(\eta_{n+\kappa})|/|\calR_n|$ and $\bar F:=\Prob[\forall n\geq 1: X_n\neq e]>0$. As $\ell_{\mathrm{supp}}$
exists, the limit $\ell_1$ also exists.
If we set $D_n:=\sum_{j=1}^{|\calR_n|}\Delta_{n+\kappa-\s{j},j}/|\calR_n|$, then (\ref{abschaetzung1}) yields the inequality \mbox{$\ell_{1}\geq
  \limsup_{n\in\mathbb{N}} D_n$}. Since the $D_n$'s are bounded, Fatou's lemma yields
$$
\ell_{1} \geq \E\Bigl[\limsup_{n\in\mathbb{N}} D_n\Bigr] \geq 
\limsup_{n\in\mathbb{N}} \E[D_n].
$$
With the help of Lemma \ref{delta-lemma1} we obtain:
\begin{eqnarray*}
\E\biggl[\frac{1}{|\calR_n|}\sum_{j=1}^{|\calR_n|}\Delta_{n+\kappa-\s{j},j} \biggr] & = &
\sum_{m=1}^{n+1} \Prob\bigl[|\calR_n|=m\bigr] \cdot
\E\biggl[\frac{1}{|\calR_n|}\sum_{j=1}^{|\calR_n|}\Delta_{n+\kappa-\s{j},j}\, \biggr|\, |\calR_n|=m \biggr]\\
& \geq & \sum_{m=1}^{n+1} \Prob\bigl[|\calR_n|=m\bigr] \cdot C = C >0.
\end{eqnarray*}
This yields
$$
\ell_{\mathrm{supp}}\geq \bar F \cdot \limsup_{n\in\mathbb{N}} \E[D_n] \geq \bar F \cdot C >0.
$$
The rest follows by (\ref{def-ell}) and $\ell\geq \ell_{\mathrm{TS}} \geq \ell_0$.
\end{proof}
We can generalize the last theorem, if we do not necessarily assume transience
of the projection $(X_n)_{n\in\mathbb{N}_0}$:
\begin{Th}
For the lamplighter random walk with respect to the infinite, finitely
generated group $G$, we have $\ell_{\mathrm{supp}}>0$ if and only if
$(X_n)_{n\in\mathbb{N}_0}$ is transient.
\end{Th}
\begin{proof}
By Theorem \ref{WoS-Th}, transience implies $\ell_{\mathrm{supp}}>0$. For the
proof of the inverse direction assume
now $\ell_{\mathrm{supp}}>0$. We have $|\mathrm{supp}(\eta_n)|\leq |B(e,R)|\cdot
|\mathcal{R}_n|$. This yields
\begin{eqnarray*}
0< \ell_{\mathrm{supp}} & = & \lim_{n\to\infty} \frac{|\mathrm{supp}(\eta_n)|}{n} \\
& \leq  & \lim_{n\to\infty} \frac{|B(e,R)|\cdot |\mathcal{R}_n|}{n} =
|B(e,R)|\cdot \Prob[\forall n\geq 1:X_n \neq e].
\end{eqnarray*}
Thus, $\Prob[\exists n\geq 1:X_n=e]<1$, that is, $(X_n)_{n\in\mathbb{N}_0}$ is transient.
\end{proof}
We get also as a consequence of Theorem \ref{WoS-Th}:
\begin{Cor}
The core rate of escape satisfies $\ell_{\mathrm{TS}}>0$. 
\end{Cor}
\begin{proof}
We have $d_{\mathrm{TS}}(\eta_n,X_n) \geq \bigl(|\mathrm{supp}(\eta_n)|-1\bigr)\cdot r_1$. Dividing both sides of the
inequality by $n$  and taking limits yields the proposed claim.
\end{proof}
Finally, we give an explicit formula for $\ell_{\mathrm{supp}}$ for the special
case of a ``Walk-Switch'' lamplighter random walk.
\par
\stepcounter{Th}
\textbf{Example \arabic{section}.\arabic{Th}:}
Suppose we are given a probability measure $\mu_0$ on $G$ with finite support
and $\langle \mathrm{supp}(\mu_0)\rangle =G$. Then a ``Walk-Switch''
lamplighter random walk over $G$ is given by the transition probabilities
$$
p\bigl( (\eta_1,x_1),(\eta_2,x_2)\bigr) :=
\frac{1}{2}\mu_0(x_1^{-1}x_2)
$$
for $x_1,x_2\in G$ and $\eta_1,\eta_2\in\mathcal{N}$ with $\eta_1=\eta_2$ or
$\eta_1\oplus \mathds{1}_{x_2}=\eta_2$. It is easy to see that
$\E[|\mathrm{supp}(\eta_n)|]=\E[|\mathcal{R}_n|]/2$, providing 
$$
\ell_{\mathrm{supp}} = \frac{1}{2} \bigl(1- \Prob[\exists n\geq 1: X_n=e]\bigr).
$$

\section{The Core Rate of Escape}

\label{SWS}
In this section we want to prove $\ell_{\mathrm{TS}}>\ell_0$ whenever $G$ is
generated as a semigroup by a symmetric set $S$ with at least three elements.
If, however, $G=\Z$ we have to make some weak assumption on the lengths of
the elements of $S$ to show $\ell_\mathrm{TS}>\ell_0$; otherwise we can
construct counterexamples where $\ell_\mathrm{TS}=\ell_0$. In this section we
may again assume $\ell_0>0$, since $\ell_\mathrm{TS}>0$.

\subsection{Groups generated by at least three elements}
\label{2ended-graphs}

In this section we assume that $S=\{s_1,\dots,s_r\}$ is symmetric with $r\geq
3$ such that there is no symmetric set
$\{s,s'\}\subseteq S$ with $G=\langle s,s'\rangle$. If $G=\langle s,s'\rangle$, then $G$ is
isomorphic to $\Z$ or \mbox{$\Z_2\ast\Z_2=\langle a,b \mid a^2=b^2=e\rangle$,}
and we treat these cases in Section \ref{z-isomorphic}. 
We want to prove that $\ell_\mathrm{\mathrm{TS}}>\ell_0$ under the
above assumption to $S$. Without loss of generality we may assume that the
elements of $S$ are ordered such that $l(s_1)\leq l(s_2)\leq \dots \leq
l(s_r)$. Our next aim is to choose three elements
$\sigma_1,\sigma_2,\sigma_3\in S$ such that $d(\sigma_k,\sigma_l)\geq \max\{
d(e,\sigma_k),d(e,\sigma_l)\}$ (with one single exception). For this purpose, we have to make a case
distinction:
\begin{enumerate}
\item[I.] If $s_1\neq s_1^{-1}$, then we define $\sigma_1:=s_1$, $\sigma_2 :=
  s_1^{-1}$ and we set $\sigma_3 := s_i$ with \mbox{$i=\min \{ k\geq 2 \mid s_k\notin \langle
s_1^{-1},s_1,s_2,\dots,s_{k-1}\rangle \}$,} that is, $\sigma_3$ is not a multiple of
$\sigma_1$ or $\sigma_2$. Note that we ensured by the above assumptions on $S$  the existence of such a $\sigma_3$.\\
\item[II.] If $s_1= s_1^{-1}$ and $s_2= s_2^{-1}$, then define $\sigma_1:=s_1$,
  $\sigma_2 :=s_2$ and $\sigma_3 := s_i$ with $i=\min \{ k\geq 3 \mid s_k\notin \langle
s_1,\dots,s_{k-1}\rangle \}$. E.g., this may happen in the case $G=\langle
  a,b,c \mid a^2=b^2=c^2=e\rangle$.\\
\item[III.] If $s_1= s_1^{-1}$ and $s_2\neq  s_2^{-1}$, then define
  $\sigma_1:=s_1$,  $\sigma_2 :=s_2$  and $\sigma_3 := s_2^{-1}$. E.g., this
  may happen in the case $G=\Z/4\Z \times \Z^3$, where $(2,0),(1,0)\in S$.
\end{enumerate}
We will see that in fact it is not relevant which one of the above cases
happens. In each of the cases we get the
following straightforward equalities, resp. inequalities:
\begin{equation}\label{bounds1}
\begin{array}{c}
d(e,\sigma_1)  =  l(\sigma_1)=r_1, \
d(e,\sigma_2)  \leq  l(\sigma_2), \ 
d(e,\sigma_3)  \leq  l(\sigma_3), \\[1ex]
d(\sigma_1,\sigma_2)  \leq  l(\sigma_1)+l(\sigma_2),\ 
d(\sigma_1,\sigma_3)  \leq  l(\sigma_1)+l(\sigma_3), \\[1ex]
d(\sigma_2,\sigma_3)  \leq  l(\sigma_2)+l(\sigma_3).
\end{array} 
\end{equation}
Moreover, we get the following (uniform) lower bounds:
\begin{Lemma} \label{bounds2}
In all the cases I, II, III,
$$
(i)\ d(e,\sigma_2)  =  l(\sigma_2),
\quad
(ii)\ d(e,\sigma_3)  =  l(\sigma_3),\quad
(iii)\ d(\sigma_1,\sigma_k)  \geq  l(\sigma_k) \textrm{ for } k\in\{2,3\}.
$$
Furthermore,
$$
(iv)\  d(\sigma_2,\sigma_3)  \geq 
\begin{cases}
l(\sigma_3), & \textrm{ in case I and II},\\
l(\sigma_1), & \textrm{ in case III}. 
\end{cases} 
$$
\end{Lemma}
\begin{proof}
Equation $(i)$ follows from
$l(\sigma_1)=l(\sigma_1^{-1})$ in the case I. In the cases II and III, we have
$s_2\neq s_1^m\in\{e,s_1\}$ for each $m\in\N$. Thus (\ref{bounds1}) yields 
equation $(i)$. In case III, equation $(ii)$ holds, since $s_2^{-1}$ is
 not a multiple of $s_1$.
For the proof of $(ii)$ in case I and II, assume $d(e,\sigma_3) < l(\sigma_3)$, that
 is, there is a path $[e,x_1,\dots,x_m=\sigma_3]$ with
 $d(x_{j-1},x_j)<l(\sigma_3)$, that is, $\sigma_3=s_i$ can be written as a
 product of elements of $s_1^{-1},s_1,\dots,s_{i-1}$, a contradiction to the minimal
 choice of $\sigma_3$. The inequalities $(iii)$ and $(iv)$ are
proved by analogous arguments. Note that the case distinction in $(iv)$ is
necessary, as the equation $s_2^2=s_1$ in case III may hold. 
\end{proof}
With the last lemma we can prove the following lemma:
\begin{Lemma}\label{lemma-distance-increase}
Let $A=\{e,\sigma_1,\sigma_2,\sigma_3\}$ and let $\varphi: \{1,2,3,4\}\to A$
be an injective function. Then in each of the cases I, II, III,
$$
d\bigl(\varphi(1),\varphi(4)\bigr) + r_1 \leq 
d\bigl(\varphi(1),\varphi(2)\bigr)+d\bigl(\varphi(2),\varphi(3)\bigr)
+d\bigl(\varphi(3),\varphi(4)\bigr).
$$
\end{Lemma}
\begin{proof}
The lemma states that for each choice of $\varphi$ a shortest
tour from $\varphi(1)$ to $\varphi(4)$ visiting $\varphi(2)$ and
$\varphi(3)$ on this trip has length at least 
\mbox{$d\bigl(\varphi(1),\varphi(4)\bigr) +l(\sigma_1)$.} Assume for the moment
that $\varphi(1)=e$ and $\varphi(4)=\sigma_3$. Then before finally reaching
$\sigma_3$ a shortest tour visiting all elements of $A$ has to pass
through $\sigma_1$ and $\sigma_2$ in this order (or first through $\sigma_2$
and then $\sigma_1$); it is not forbidden to visit $e$ or $\sigma_3$ twice. This tour has a length of at least
$$
d(e,\sigma_1)+d(\sigma_1,\sigma_2) +d(\sigma_2,\sigma_3) \geq
\begin{cases}
l(\sigma_1) + l(\sigma_2) + l(\sigma_3), & \textrm{in cases I and II},\\
2\,l(\sigma_1)+ l(\sigma_2), & \textrm{in case III}.
\end{cases}
$$
But $d(e,\sigma_3)\leq l(\sigma_3)$ in case I and II, and $d(e,\sigma_3)\leq
l(\sigma_2)$ in case III. Thus, the claim follows for the specific choice of
$\varphi$ with $\varphi(1)=e, \varphi(2)=\sigma_1,\varphi(3)=\sigma_2,\varphi(4)=\sigma_3$. For all
other choices of $\varphi$ the same result follows; compare with Figure
\ref{distance-compare}: the first four columns build a case distinction for the
choice of $\varphi$ (we use symmetries!), the fifth column gives an upper bound for
$d\bigl(\varphi(1),\varphi(4)\bigr)$; the sixth column gives a lower bound for the right
hand side of the inequality in the lemma
and the last column is a lower bound for the difference between the sixth and
fifth column $\bigl($recall that $l(\sigma_1)\leq l(\sigma_2) \leq
l(\sigma_3)\bigr)$. For case III, we only summarize the different possibilities,
where $\varphi(i-1)=\sigma_2,\, \varphi(i)=\sigma_3$ or
$\varphi(i-1)=\sigma_3,\, \varphi(i)=\sigma_2$, where we use Lemma \ref{bounds2}(iv); the lower and upper bounds for
all other choices for $\varphi$ conincide with case I and II. Compare also with
Figure \ref{triangles}, where the labels on the dotted lines are the lower
bounds for the distances between two points.
This proves the lemma.
\begin{figure}[htp]
Cases I and II:\\[1ex]
\begin{tabular}{c|c|c|c|c|c|c}
$\varphi(1)$ & $\varphi(2)$ & $\varphi(3)$ & $\varphi(4)$ & 
$d\bigl(\varphi(1),\varphi(4)\bigr) \leq$ & Right Side $\geq$ & Difference $\geq$ \\
\hline
$e$ & $\sigma_1$ & $\sigma_2$ & $\sigma_3$ & $l(\sigma_3)$ &
$l(\sigma_1)+l(\sigma_2)+l(\sigma_3)$ & $l(\sigma_1)+l(\sigma_2)$\\
$e$ & $\sigma_1$ & $\sigma_3$ & $\sigma_2$ & $l(\sigma_2)$ &
$l(\sigma_1)+2\,l(\sigma_3)$ & $l(\sigma_1)+l(\sigma_3)$\\
$e$ & $\sigma_2$ & $\sigma_1$ & $\sigma_3$ & $l(\sigma_3)$ &
$2\,l(\sigma_2)+l(\sigma_3)$ & $2\,l(\sigma_2)$\\
$e$ & $\sigma_3$ & $\sigma_1$ & $\sigma_2$ & $l(\sigma_2)$ &
$l(\sigma_2)+2\,l(\sigma_3)$ & $2\,l(\sigma_3)$\\
$e$ & $\sigma_2$ & $\sigma_3$ & $\sigma_1$ & $l(\sigma_1)$ &
$l(\sigma_2)+2\,l(\sigma_3)$ & $2\,l(\sigma_3)$\\
$e$ & $\sigma_3$ & $\sigma_2$ & $\sigma_1$ & $l(\sigma_1)$ &
$l(\sigma_2)+2\,l(\sigma_3)$ & $2\,l(\sigma_3)$\\
\hline
$\sigma_1$ & $e$ & $\sigma_2$ & $\sigma_3$ & $l(\sigma_1)+l(\sigma_3)$ &
$l(\sigma_1)+l(\sigma_2)+l(\sigma_3)$ & $l(\sigma_2)$\\
$\sigma_1$ & $e$ & $\sigma_3$ & $\sigma_2$ & $l(\sigma_1)+l(\sigma_2)$ &
$l(\sigma_1)+2\,l(\sigma_3)$ & $l(\sigma_3)$\\
$\sigma_1$ & $\sigma_2$ & $e$ & $\sigma_3$ & $l(\sigma_1)+l(\sigma_3)$ &
$2\,l(\sigma_2)+l(\sigma_3)$ & $l(\sigma_2)$\\
$\sigma_1$ & $\sigma_3$ & $e$ & $\sigma_2$ & $l(\sigma_1)+l(\sigma_2)$ &
$l(\sigma_2)+2\,l(\sigma_3)$ & $l(\sigma_3)$\\
\hline
$\sigma_2$ & $e$ & $\sigma_1$ & $\sigma_3$ & $l(\sigma_2)+l(\sigma_3)$ &
$l(\sigma_1)+l(\sigma_2)+l(\sigma_3)$ & $l(\sigma_1)$ \\
$\sigma_2$ & $\sigma_1$ & $e$ & $\sigma_3$ & $l(\sigma_2)+l(\sigma_3)$ &
$l(\sigma_1)+l(\sigma_2)+l(\sigma_3)$ & $l(\sigma_1)$ 
\end{tabular}\\[1ex]
Case III:\\[1ex]
\begin{tabular}{c|c|c|c|c|c|c}
$\varphi(1)$ & $\varphi(2)$ & $\varphi(3)$ & $\varphi(4)$ & 
$d\bigl(\varphi(1),\varphi(4)\bigr) \leq$ & Right Side $\geq$ & Difference $\geq$ \\
\hline
$e$ & $\sigma_1$ & $\sigma_2$ & $\sigma_3$ & $l(\sigma_2)$ &
$2\,l(\sigma_1)+l(\sigma_2)$ & $2\,l(\sigma_1)$\\
$e$ & $\sigma_1$ & $\sigma_3$ & $\sigma_2$ & $l(\sigma_2)$ &
$2\,l(\sigma_1)+l(\sigma_2)$ & $2\,l(\sigma_1)$\\
$e$ & $\sigma_2$ & $\sigma_3$ & $\sigma_1$ & $l(\sigma_1)$ &
$l(\sigma_1)+2\,l(\sigma_2)$ & $2\,l(\sigma_2)$\\
$e$ & $\sigma_3$ & $\sigma_2$ & $\sigma_1$ & $l(\sigma_1)$ &
$l(\sigma_1)+2\,l(\sigma_2)$ & $2\,l(\sigma_2)$\\
\hline
$\sigma_1$ & $e$ & $\sigma_2$ & $\sigma_3$ & $l(\sigma_1)+l(\sigma_2)$ &
$2\,l(\sigma_1)+l(\sigma_2)$ & $l(\sigma_1)$\\
$\sigma_1$ & $e$ & $\sigma_3$ & $\sigma_2$ & $l(\sigma_1)+l(\sigma_2)$ &
$2\,l(\sigma_1)+l(\sigma_2)$ & $l(\sigma_1)$\\
\end{tabular}
\caption{Comparision for the choice of $\varphi$}
\label{distance-compare}
\end{figure}
\begin{figure}[htp]
\includegraphics[width=3cm]{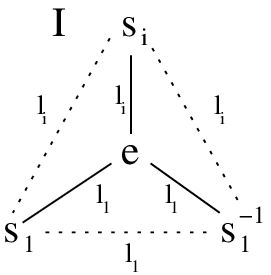}\hspace{5mm}
\includegraphics[width=3cm]{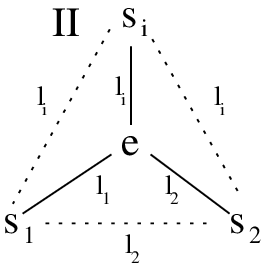}\hspace{5mm}
\includegraphics[width=3cm]{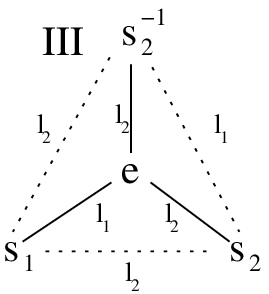}\hspace{5mm}
\caption{Distances between $e$, $\sigma_1$, $\sigma_2$, $\sigma_3$ with
  $l_1=l(s_1)$, \mbox{$l_2=l(s_2)$} and $l_i=l(s_i)$.}
\label{triangles}
\end{figure}
\end{proof}
In other words, the lemma states that a shortest tour starting at some $a\in A$
 visiting all other elements of $A$ and finishing at some $a'\in A$ has length of at least
$d(a,a')+r_1$. We will now apply this lemma independently of which of the cases I,
II, III applies. 
For $y\in G$ let $B_y:=\{y,y\sigma_1,y\sigma_2,y\sigma_3\}$ and let be $x_y\in
G\setminus B_y$. Obviously, for each choice of $w,z\in B_y$,
\begin{equation}\label{distance-increase1}
d(e,x_y) \leq d(e,w) + d(w,z) + d(z,x_y).
\end{equation}
Let $\mathcal{F}$ be the set of all injective functions $\varphi:\{1,2,3,4\}\to
B_y$. Then the last lemma and (\ref{distance-increase1}) yield the following inequality:
\begin{eqnarray}
&&d_\mathrm{TS}\bigl(\mathds{1}_y \oplus \mathds{1}_{y\sigma_1} \oplus
\mathds{1}_{y\sigma_2} \oplus \mathds{1}_{y\sigma_3},x\bigr) \nonumber \\[2ex]
&\geq &\min_{\varphi\in\mathcal{F}} 
\biggl\lbrace d\bigl(e,\varphi(1)\bigr) + \sum_{i=1}^3 d\bigl(\varphi(i),\varphi(i+1)\bigr)+d\bigl(\varphi(4),x\bigr)\biggr\rbrace \nonumber\\[2ex]
& \geq & d(e,x) + r_1.\label{distance-increase2}
\end{eqnarray}
We now come back to our lamplighter random walk. Our next aim is to bound $d_\mathrm{TS}(Z_n)$ from
below with the help of (\ref{distance-increase2}), independently of which of
the cases I, II, III holds. For this
purpose, define \textit{hitting times}
$$
\t{1}:=0, \, \textrm{and } \t{k}:=\min \biggl\lbrace m\in \N \,\biggl|\, m>\t{k-1}, X_m\notin \bigcup_{j=1}^{k-1} B\bigl(X_{\t{j}},2l(\sigma_3)\bigr)\biggr\rbrace \ \textrm{for } k\geq 2,
$$
that is, $\t{k}$ is the first instant of time after time $\t{k-1}$ for which the random walk leaves the finite set $\bigcup_{j=0}^{k-1} B\bigl(X_{\t{j}},2l(\sigma_3)\bigr)$.
By transience, $\t{k}<\infty$ almost surely. Furthermore, we write $\H{k}:=X_{\t{k}}$ and $\calR_n':=\{X_{\t{j}}\,|\, j\in\N \textrm{ with } \t{j}\leq n\}$.
Observe that
\begin{equation}\label{empty-set}
\{\H{k},\H{k}\sigma_1,\H{k}\sigma_2,\H{k}\sigma_3\}\cap
\{\H{l},\H{l}\sigma_1,\H{l}\sigma_2,\H{l}\sigma_3\} =\varnothing
\end{equation}
for $k\neq l$.  
The idea is to investigate, if enough lamps are on in $B\bigl(X_{\t{j}},l(\sigma_3)\bigr)$
such that we have
$d_\mathrm{TS}(\eta,y)> d(e,y)$, where $y\in G$ and 
$\mathrm{supp}(\eta)$ is a subset of this ball. Our aim is
to construct deviations to establish such a situation for each of these balls with a
strict positive probability independently \mbox{of $k$.} See
Figure \ref{exittimes-figure1}.
\begin{figure}[thb]
\includegraphics[width=4cm]{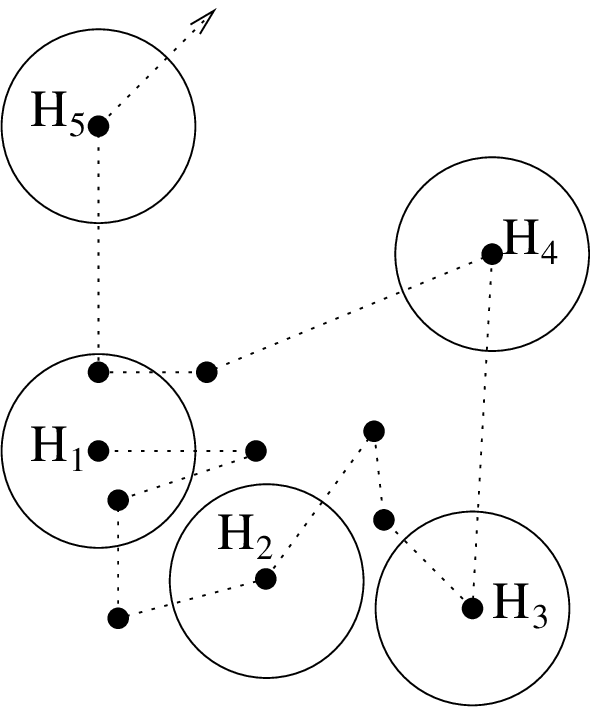}
\caption{Hitting points}
\label{exittimes-figure1}
\end{figure}
\par
For $n\in\N_0$, $k\in \N$, define
$$
\Delta_{n,k} := 
\begin{cases}
r_1, & \textrm{if }
\{\H{k},\H{k}\sigma_1,\H{k}\sigma_2,\H{k}\sigma_3\} \subseteq
\mathrm{supp}(\eta_{\t{k}+n}),\\
0, & \textrm{otherwise}.
\end{cases}
$$
If $n\geq \t{k}$ and $\Delta_{n-\t{k},k}=1$, a shortest tour
from $e$ to $X_n$ visiting each element of $\mathrm{supp}(\eta_{n})$ has
to visit in particular each element of $\bigl\lbrace
\H{k},\H{k}\sigma_1,\H{k}\sigma_2,\H{k}\sigma_3\bigr\rbrace$. But by Lemma
\ref{lemma-distance-increase} and (\ref{distance-increase2}) this means
that
$$
d_\mathrm{TS}(Z_n) \geq d(e,X_n) + r_1.
$$
Due to (\ref{empty-set}), iterated applications of
the triangular inequality and \mbox{Lemma \ref{lemma-distance-increase}} yield
\begin{equation}\label{dTS-estimate}
d_\mathrm{TS}(Z_{n}) \geq d(e,X_{n})
+ \sum_{j=1}^{|\calR_n'|} \Delta_{n-\t{j},j}.
\end{equation}
Our next aim is to bound $\Prob[\Delta_{n,k}=r_1]$ for $n$ big enough, and thus
$\E[\Delta_{n,k}]$, uniformly from below by a non-zero constant. 
\begin{Prop}\label{delta2}
There exist $\lambda\in\N$ and $D>0$ such that $\E[\Delta_{n,k}] \geq D$ for all $k,n\in\mathbb{N}$
with $n\geq \lambda$.
\end{Prop}
\begin{proof}
By transience, there is $y\in G$ with $d(e,y)\geq R+l(\sigma_3)$ such that
$$
\hat p := \Prob_{(\eta',y)}\bigl[\forall n\geq 1:X_n\notin
B\bigl(e,R+l(\sigma_3)\bigl)\bigr]>0 \quad \textrm{ for all } \eta'\in\calN.
$$
For each subset $A\subseteq \{e,\sigma_1,\sigma_2,\sigma_3\}$ there is
$\lambda_A\in\N$ such that 
$$
D_A:=\mu^{(\lambda_A)}\Bigl(\Bigl(\sum_{w\in A}
\mathds{1}_w,y\Bigr)\Bigr)>0.
$$
We sum over all possibilities for the hitting point $\H{k}$. Assume now for a moment that $x$ is
the hitting point. Thus, the probability of walking from $x$ with configuration $\eta$ to $xy$ such that the lamps at $x,x\sigma_1,x\sigma_2,x\sigma_3$ are on when the lamplighter
 reaches $xy$ is at least $D':=\min\bigl\lbrace D_A | A\subseteq
 \{e,\sigma_1,\sigma_2,\sigma_3\}\bigr\rbrace$. The lamplighter will not return to
 $B\bigl(x,R+l(\sigma_3)\bigr)$ when starting at $xy$ with a positive
 probability, namely with a probability of at least $\hat p$. We now obtain for $n\geq \lambda:=\max\bigl\lbrace \lambda_A | A\subseteq
 \{e,\sigma_1,\sigma_2,\sigma_3\}\bigr\rbrace$:
\begin{eqnarray*}
\E[\Delta_{n,k}] & = & r_1\cdot \Prob[\Delta_{n,k}=r_1] \\[1ex]
&= & r_1\cdot \sum_{(\eta,x)\in \calL} \sum_{m\geq 0}
\Prob\bigl[\t{k}=m,Z_{m}=(\eta,x),
\{x,x\sigma_1,x\sigma_2,x\sigma_3\}\subseteq 
\mathrm{supp}(\eta_{m+n})\bigr]\\
&\geq & r_1\cdot \sum_{(\eta,x)\in \calL} \sum_{m\geq 0}
\Prob\bigl[\t{k}=m,Z_{m}=(\eta,x)\bigr] \cdot D' \cdot \hat p \\
& = & r_1 \cdot D' \cdot \hat p =: D. 
\end{eqnarray*}
\end{proof}
Now we can summarize:
\begin{Th}\label{dTS-th}
For the lamplighter random walk on the infinite, finitely generated group
$G$, generated as a semigroup by the symmetric set $S$ such that $G\neq \langle
s_k,s_l\rangle$ for $s_k,s_l\in S$, 
$$
\ell \geq \ell_{\mathrm{TS}} > \ell_0.
$$
\end{Th}
\begin{proof}
In view of inequality
(\ref{dTS-estimate}), Fatou's lemma and Proposition \ref{delta2} give
\begin{eqnarray*}
\ell_{\mathrm{TS}} & = & \int \lim_{n\to\infty}
\frac{d_{\mathrm{TS}}\bigl(Z_{n+\lambda}\bigr)}{n}\,d\Prob \\
&\geq & \int \limsup_{n\in\N} \biggl( \frac{d(e,X_{n+\lambda})}{n}+ \frac{1}{n}
\sum_{j=1}^{|\calR_n'|} \Delta_{n+\lambda-\t{j},j}\biggr)\,d\Prob \\
&\geq & \ell_0 + \limsup_{n\in\N} \int \biggl( \frac{1}{n}
\sum_{j=1}^{|\calR_n'|} \Delta_{n+\lambda-\t{j},j}\biggr)\,d\Prob \\
&=& \ell_0 + \limsup_{n\in\N} \sum_{k=1}^{n+1} \Prob\bigl[|\calR_n'|=k\bigr] \cdot
\E\biggl[\frac{1}{n} \sum_{j=1}^{|\calR_n'|} \Delta_{n+\lambda-\t{j},j}\,\biggl|\,
|\calR_n'|=k\biggr].
\end{eqnarray*}
This provides
$$
\ell_{\mathrm{TS}} \geq
 \ell_0 + \limsup_{n\in\N} \sum_{k=1}^{n+1} \Prob\bigl[|\calR_n'|=k\bigr] \cdot
\frac{k}{n}\cdot D
=  \ell_0 + D \cdot \limsup_{n\in\N} \frac{\E[|\calR_n'|]}{n}.
$$
By Kingman's subadditive ergodic theorem $\E[|\calR_n'|]/n$ converges, and due
to the inequality $|\calR_n|\leq |B\bigl(e,2l(\sigma_3)\bigr)|\cdot |\calR_n'|$
its limit is bounded from below by some constant $\bar D>0$, completing the 
proof.
\end{proof}

\subsection{$\Z$-isomorphic Groups}

\label{z-isomorphic}
In this section we consider the remaining case, where $G$ is generated as a
semigroup by two
elements $s_k,s_l\in S$, that is, $G\simeq \Z$ or
$G\simeq \Z_2\ast\Z_2=\langle a,b \mid a^2=b^2=e\rangle$. For the sake of
completeness we prove the following lemma:
\begin{Lemma}
Any irreducible random walk $(X_n)_{n\in\N_0}$ on $G=\Z_2\ast\Z_2$ governed by
a probability measure $\mu_0$ with finite first moment
is recurrent. In particular, $\ell_{\mathrm{TS}}=0$.
\end{Lemma}
\begin{proof}
Observe that
$$
Z := \bigl\lbrace (ab)^z \mid z\in\Z \bigr\rbrace \simeq \Z
$$
is a subgroup of $G$, which has index $2$. We identify from now on the elements
of $Z$ with integers and write $\overline{Z}$ for the
complement $G\setminus Z$. Consider the stopping times $T_0:=0$,
$T_n:=\min\{m\in\N \mid X_m\in Z, m>T_{m-1}\}$ for $m\in\N$, and define $Y_n:=X_{T_n}$. Then
the random walk $(X_n)_{n\in\N_0}$ on $G$ is recurrent if and only if $(Y_n)_{n\in\N_0}$ is
recurrent. We now compute the expectation of the drift $D_n:=Y_{n}-Y_{n-1}$, which is independent of $n$. We write
$\mu_Z(z)=\mu_0(z)/\mu_0(Z)$ for every $z\in \Z$ and 
$d_Z=\sum_{w\in Z} w\,\mu_Z(w)$.
Observe that for each $w,w'\in\overline{Z}$ we have $ww'\in Z$ and $ww'=-w'w$.
Now we can compute the expected drift by distinguishing how one can walk from
$Y_{n-1}$ to $Y_n$:
\begin{eqnarray*}
\E[D_n] & = & \sum_{w\in Z} w\,\mu_0(w) +
\sum_{n\geq 1}
\sum_{w_0,w_n\in \overline{Z}}  \sum_{w_1,\dots,w_{n-1}\in Z}
(w_0-w_1-\dots - w_{n-1}+ w_n) \cdot\\
&& \quad \hspace{6.5cm} \cdot \mu_0(w_0) \mu_0(w_1) \cdots \mu_0(w_{n-1})
\mu_0(w_n)\\[1ex]
&=& \mu_0(Z)\cdot d_Z + \sum_{n\geq 1} \mu_0(Z)^{n-1} \Bigl(
 \sum_{w_0,w_n\in \overline{Z}} (w_0+w_n)\mu_0(w_0)\mu_0(w_n) +\\
&& \quad \quad
 +\bigl(1-\mu_0(Z)\bigr)^2 \cdot \sum_{w_1,\dots,w_{n-1}\in Z} (-w_1-\dots -w_{n-1})
 \mu_Z(w_1)\dots \mu_Z(w_{n-1})\Bigr)\\
&=& \mu_0(Z)\cdot d_Z - \bigl(1-\mu_0(Z)\bigr)\cdot d_Z \sum_{n\geq 1}
(n-1)\,\mu_0(Z)^{n-1}\, \bigl(1-\mu_0(Z)\bigr)\\
&=& \mu_0(Z)\cdot d_Z - \bigl(1-\mu_0(Z)\bigr)\cdot d_Z \cdot \mu_0(Z)/\bigl(1-\mu_0(Z)\bigr) = 0.
\end{eqnarray*}
By Chung-Ornstein \cite{chung-ornstein} it follows that
$(Y_n)_{n\in\N_0}$ is recurrent, and thus $(X_n)_{n\in\N_0}$ is recurrent. From Dyubina
\cite{dyubina} follows $\ell_{\mathrm{TS}}=0$. 
\end{proof}
We want to show that \mbox{Theorem \ref{dTS-th}} holds also in the case $G=\Z$ under suitable assumptions on the lengths
of the elements of $S$. Note that we may assume again that $\ell_0>0$ and that
$|S|\geq 3$. If $|S|=2$,
then $G$ is isomorphic to $\mathbb{Z}=\langle S\rangle$ with $S=\{-1,1\}$,
whose Cayley graph is the infinite line, on which one can easily show that the lamplighter random walk has the same speed
as its projection \mbox{onto $G$}; see also Bertacchi \cite{bertacchi}. Thus, we
only have to take a closer look on $\mathbb{Z}$ generated by a symmetric set
$S$ with $-1,1\in S$ and $|S|\geq 3$. Observe that if $\pm 1\notin S$, then we
may apply the results of the previous section.
 Furthermore, assume that there is $s\in S\setminus
\{\pm 1\}$ with $l(s)<|s|\cdot l(1)$; otherwise we are more or less in the
situation of $S=\{\pm 1\}$, see the end of this section. Moreover, we may assume that $d(0,1)=d(0,-1)=l(1)$; otherwise
$[x_0=0,x_1=1]$ is not a shortest path from $0$ to $1$, that is, $S\setminus\{-1,1\}$
provides the same metric as the metric induced by $S$, that is, we may apply
the results of Section \ref{2ended-graphs} in such case. Due to the same
argument we may assume that the only shortest path from $0$ to $1$ is $[x_0=0,x_1=1]$.
\par
We proceed similarily to Section \ref{2ended-graphs}. We make a case
distinction and define:
\begin{itemize}
\item[I.]
If there is $s\in S\setminus \{ \pm 1\}$ such that $r_1=l(s)<l(1)$, then define $\sigma_1=s$, $\sigma_2=s^{-1}$
and $\sigma_3=1$. \\
\item[II.] Otherwise we set $\sigma_1=1$, $\sigma_2=-1$ and $\sigma_3=s$, where $s\in\N\cap (S\setminus \{1\})$ such that $l(s)<|s|\cdot l(1)$ and $l(s)\leq l(s')$ for all $s'\in S$ with $l(s')<|s'|\cdot l(1)$.
\end{itemize}
In case I we have trivially $d(0,s)=d(0,s^{-1})=l(s)=r_1$ and $d(s,s^{-1})\geq
l(s)=r_1$, while in case II we have $d(0,1)=d(0,-1)=l(1)=r_1$ and $d(1,-1)\geq l(1)$. Moreover:
\begin{Lemma} 
We have the following equations and lower bounds:
\begin{enumerate}
\item In case I there is some $\varepsilon_0 >0$ such that
\begin{itemize}
\item[$(i)$] $d(0,1)=l(1)$,
\item[$(ii)$] $d(s,1)\geq l(1)-l(s)+\varepsilon_0$,
\item[$(iii)$] $d(s^{-1},1)\geq l(1)-l(s)+\varepsilon_0.$
\end{itemize}
\item In case II there is some $\varepsilon_0 >0$ such that
\begin{itemize}
\item[$(i)$] $d(0,s)=l(s)$,
\item[$(ii)$] $d(1,s)\geq l(s)-l(1)+\varepsilon_0$,
\item[$(iii)$] $d(-1,s)\geq l(s)-l(1)+\varepsilon_0.$
\end{itemize}
\end{enumerate}
\end{Lemma}
\begin{proof}
Equation (1).$(i)$ holds by the assumption made above.
For the proof of $(1).(ii)$ assume that $d(s,1)\leq l(1)-l(s)$. Then there is a
shortest path \mbox{$[s,x_1,\dots,x_n=1]$} with $d(x_{i-1},x_i)\leq l(1)-l(s)$, that is,
$x_{i-1}^{-1}x_i\neq \pm 1$. But this means that there is another path
from $0$ to $1$ of length at most $l(1)$, namely $[0,s,x_1,\dots,x_n]$, distinct from $[0,1]$, a contradiction to the 
assumptions above. As $S$ is finite, existence of $\varepsilon_0$ is ensured. Inequality
$(1).(iii)$ is proved analogously.
\par
Assume that equation $(2).(i)$ does not hold. This implies that there is a
shortest path \mbox{$[0,x_1,\dots,x_n=s]$} from $0$ to $s$ with $d(x_{i-1},x_i)<l(s)$, that is, we may
assume $x_{i-1}^{-1}x_i= \pm 1$ by minimality of $l(s)$. But this implies
$l(s)>d(0,s)=|s|\cdot l(1)$, a contradiction to the choice of $s$. To prove
$(2).(ii)$ assume $d(1,s)\leq l(s)-l(1)$, that is, there is a
shortest path $[1,x_1,\dots,x_n=s]$ with $d(x_{i-1},x_i)\leq l(s)-l(1)$, that is,
we may assume $x_{i-1}^{-1}x_i= \pm 1$ by minimality of $l(s)$. But this provides now
$l(s)\geq l(1)+d(1,s) = |s|\cdot l(1)$, a contradiction to the choice of $s$.
Inequality $(2).(iii)$ is proved analogously.
\end{proof}
We get the analogue to Lemma \ref{lemma-distance-increase}:
\begin{Lemma}\label{lemma-distance-increase2}
Let be $A=\{0,\sigma_1,\sigma_2,\sigma_3\}$ and let $\varphi: \{1,2,3,4\}\to A$
be an injective function. Then in each of the cases I and II,
$$
d\bigl(\varphi(1),\varphi(4)\bigr) + \min\{\varepsilon_0,l(s),l(1)\} \leq 
d\bigl(\varphi(1),\varphi(2)\bigr)+d\bigl(\varphi(2),\varphi(3)\bigr)
+d\bigl(\varphi(3),\varphi(4)\bigr).
$$
\end{Lemma}
\begin{proof}
The proof works analogously to the proof of Lemma
\ref{lemma-distance-increase}; compare with \mbox{Figure \ref{distance-compare-Z}} for
the comparision of the distances in case I. The inequality for case II follows
analogously by symmetry.
\begin{figure}[htp]
{\small
\begin{tabular}{c|c|c|c|c|c|c}
$\varphi(1)$ & $\varphi(2)$ & $\varphi(3)$ & $\varphi(4)$ & 
$d\bigl(\varphi(1),\varphi(4)\bigr) \leq$ & Right Side $\geq$ & Difference $\geq$ \\
\hline
$0$ & $s$ & $s^{-1}$ & $1$ & $l(1)$ &
$2\,l(s)+(l(1)-l(s)+\varepsilon_0)$ & $l(s)+\varepsilon_0$\\
$0$ & $s^{-1}$ & $s$ & $1$ & $l(1)$ &
$2\,l(s)+(l(1)-l(s)+\varepsilon_0)$ & $l(s)+\varepsilon_0$\\
$0$ & $1$ & $s^{-1}$ & $s$ & $l(s)$ &
$l(1)+(l(1)-l(s)+\varepsilon_0)+l(s)$ & $l(1)+\varepsilon_0$\\
$0$ & $s^{-1}$ & $1$ & $s$ & $l(s)$ &
$l(s)+2(l(1)-l(s)+\varepsilon_0)$ & $\varepsilon_0$\\
$0$ & $1$ & $s$ & $s^{-1}$ & $l(s)$ &
$l(1)+(l(1)-l(s)+\varepsilon_0)+l(s)$ & $l(1)+\varepsilon_0$\\
$0$ & $s$ & $1$ & $s^{-1}$ & $l(s)$ &
$l(s)+2(l(1)-l(s)+\varepsilon_0)$ & $2\varepsilon_0$\\
\hline
$s$ & $0$ & $1$ & $s^{-1}$ & $2\,l(s)$ &
$l(s)+l(1)+(l(1)-l(s)+\varepsilon_0)$ & $\varepsilon_0$\\
$s$ & $1$ & $0$ & $s^{-1}$ & $2\,l(s)$ &
$(l(1)-l(s)+\varepsilon_0)+l(1)+l(s)$ & $\varepsilon_0$\\
$s$ & $0$ & $s^{-1}$ & $1$ & $l(s)+l(1)$ &
$2\,l(s)+(l(1)-l(s)+\varepsilon_0)$ & $\varepsilon_0$\\
$s$ & $s^{-1}$ & $0$ & $1$ & $l(s)+l(1)$ &
$2l(s)+l(1)$ & $l(s)$\\
\hline
$s^{-1}$ & $0$ & $s$ & $1$ & $l(s)+l(1)$ &
$2\,l(s)+(l(1)-l(s)+\varepsilon_0)$ & $\varepsilon_0$ \\
$s^{-1}$ & $s$ & $0$ & $1$ & $l(s)+l(1)$ &
$2\,l(s)+l(1)$ & $l(s)$ 
\end{tabular}}\\[1ex]
\caption{Comparision for the choices of $\varphi$ in case I}
\label{distance-compare-Z}
\end{figure}
\end{proof}
Now we can conlude:
\begin{Cor}
For the lamplighter random walk on $G=\Z$, generated as a semigroup by the
symmetric set $S$ such that $-1,1\in S$, $|S|\geq 3$ and $l(s)<|s|\cdot l(1)$ for some
$s\in S\setminus\{-1,1\}$,
$$
\ell \geq \ell_{\mathrm{TS}}> \ell_0.
$$
\end{Cor}
\begin{proof}
Due to Lemma \ref{lemma-distance-increase2} the proof follows analogously to
the considerations of Section \ref{2ended-graphs}, where we redefine $\Delta_{n,k}$
by
$$
\Delta_{n,k} := 
\begin{cases}
\min\{\varepsilon_0,l(s),l(1)\}, & \textrm{if }
\{\H{k},\H{k}\sigma_1,\H{k}\sigma_2,\H{k}\sigma_3\} \subseteq
\mathrm{supp}(\eta_{\t{k}+n}),\\
0, & \textrm{otherwise}.
\end{cases}
$$
\end{proof}
We now explain the necessity of having some $s\in S\setminus\{-1,1\}$ with
$l(s)<|s|\cdot l(1)$. If this assumption is not satisfied, then the metric on
$G=\Z$ is $d(x,y)=r_1 \cdot |x-y|$, that is, we have the natural metric on
$\Z$ if $r_1=1$. In this case the lamplighter random walk
has the same speed as its projection onto the group $G$. E.g., consider $G=\mathbb{Z}$
generated by $S=\bigl\lbrace \pm 1, \pm 2, \pm 3\bigr\rbrace$ with $l(\pm 1)=1$,
$l(\pm 2)=3$ and $l(\pm 3)=5$.
 Observe that
$[0,1,2,\dots,z]$ for $z>0$ is a shortest path from $0$ to
$z$. Let be $p\in(1/2;1)$. We equip 
$\Z_2 \wr \Z$ with a transient random walk defined by the following transition
probabilities:
\par
\begin{center}
$\mu(\mathbf{0},1)  =  \mu(\mathds{1}_0,1) =\mu(\mathbf{0},2) =\mu(\mathds{1}_0,2)
= \mu(\mathbf{0},3) =\mu(\mathds{1}_0,3) = p/6$,\\[1ex]
$\mu(\mathbf{0},-1) =  \mu(\mathds{1}_0,-1) =\mu(\mathbf{0},-2) =\mu(\mathds{1}_0,-2)
= \mu(\mathbf{0},-3) =\mu(\mathds{1}_0,-3) = (1-p)/6$.
\end{center}
\par
Thus, $d(e,X_n)/n$ converges almost surely to $2p-1$.  Analogously to the case $G=\mathbb{Z}=\langle \pm 1\rangle$, it can be
shown that the lamplighter does not escape faster than its projection on $\Z$,
that is, we have $\ell_{\mathrm{TS}}=\ell_0$.

\section{Remarks}

\label{remarks}

\subsection{Generalization to Transitive Graphs and Markovian Distance}

The results of Sections \ref{WoS} and \ref{SWS} can be generalized to 
transient
lamplighter random walks on transitive, connected, locally finite graphs, which
are not necessarily Cayley graphs of finitely generated groups. Again, it is
assumed that the lamplighter random walk's projection onto the base graph is
transient. The results of the previous sections also apply in this case, if
graph automorphisms leave the
lamplighter random walk operator invariant; compare with
Gilch \cite{gilch}.
\par
One can also investigate the rate of escape with respect to the
\textit{Markovian distance} $d_\Prob\bigl((\eta,x),(\eta',x')\bigr)$ on $\Z_2
\wr G$, which is given by 
$$
\min 
\left\lbrace \sum_{i=1}^n d(x_{i-1},x_i)\, \biggl| \, 
\begin{array}{c}
n\in\N_0, \textrm{there are } x_0,x_1,\dots x_{n}\in G \textrm{ with }
x_0=x,x_n=x'   \\
\textrm{such that }
\Prob_{(\eta,x)}\bigl[X_1=x_1,\dots,X_n=x',\eta_n=\eta'\bigr]>0
\end{array}
\right\rbrace.
$$
The limit $\ell_\Prob:=\lim_{n\to\infty}
d_\Prob\bigl((\mathbf{0},e),Z_n\bigr)/n$ exists almost surely by Kingman's
subadditive ergodic theorem and is almost surely constant. It can be shown that
if the Cayley graph of $G$ has infinitely many ends, then with respect to the Markovian distance the lamplighter escapes faster to
infinity than its projection onto $G$. If the
assumption is dropped one can find counterexamples such that the
lamplighter is not faster; e.g. if $G=\Z\times \Z_2$.

\subsection{Multi-State Lamps}

The presented techniques for proving the acceleration of the lamplighter random
walks can also be applied to the case that there are more possible lamp
states encoded by elements of $\Z/r\Z$ with $r>2$. In this case one may assign
lengths to a set of generators of $\Z/r\Z$. 
Then the presented results can be proved analogously.

\subsection{Greenian Distance}

Another metric on $G$ is given by the \textit{Greenian distance}
$$
d_{\mathrm{Green}}(x,y):= -\ln \Prob_x[T_y<\infty],
$$
where $T_y$ is the hitting time of $y\in G$. Analogously, we can define the Greenian metric for the random walk on $\Z_2\wr G$. These metrics are
not path metrics induced by lengths on the set of generators. Benjamini and
Peres \cite{benjamini-peres94} proved that the entropy and the rate of escape
w.r.t. the Greenian distance of random walks on finitely generated groups with
finite support are equal.
Blach\`ere, Ha\"issinsky and
Mathieu \cite{blachere-haissinsky-mathieu} generalized this result to random
walks on countable groups. If the random walk on $G$ is
governed by a probability measure $\mu_0$ with $\langle \mathrm{supp}\,
\mu_0\rangle =G$, then
the entropy of the
lamplighter random walk on $\Z_2\wr G$ is strictly bigger than the entropy of the random walk's
projection onto $G$, because the Poisson boundary of the lamplighter random walk
projects non-trivially onto the one of the random walk on the base graph;
compare with Kaimanovich and Vershik \cite[Theorem 3.2]{kaimanovich-vershik}.
It follows that with respect to the
Greenian distance the lamplighter random walk is faster than its projection
\mbox{onto $G$.}

\section*{Acknowledgements}

The author is grateful to Wolfgang Woess and Sebastian M\"uller for numerous
discussions, and also to the anonymous referee, who gave several useful remarks
regarding content.

\bibliographystyle{abbrv}
\bibliography{literatur}

\end{document}